\documentclass[journal,english]{IEEEtran}
\usepackage{times}
\usepackage{amsmath}
\usepackage{graphicx}
\usepackage{amssymb}
\usepackage{float}
\usepackage{cite}

\providecommand{\tabularnewline}{\\}
\floatstyle{ruled}
\newfloat{algorithm}{tbp}{loa}
\floatname{algorithm}{Algorithm}

\usepackage{algorithmic}

\newcommand{\fftfunc}[3]{\operatorname{\mathbf{#1}}_{#2}\left(#3\right)}

\newcommand{\fftfuncl}[4]
{\operatorname{\mathbf{#1}}_{#2}^{#3}\left(#4\right)}

\newcommand{\fftfunclT}[4]
{\operatorname{\mathbf{#1}}_{#2}^{(\mathrm{T})\,#3}\left(#4\right)}

\newcommand{\splitfft}[2]{\fftfunc{splitfft}{#1}{#2}}
\newcommand{\rsplitfft}[2]{\fftfunc{rsplitfft}{#1}{#2}}

\newcommand{\nsplitfftds}[2]{\fftfunc{newfftS}{#1}{#2}}

\newcommand{\nsplitfftdsfour}[2]{\fftfunc{newfftS4}{#1}{#2}}

\newcommand{\rnsplitfftdsfour}[2]{\fftfunc{rnewfftS4}{#1}{#2}}

\newcommand{\nsplitfftdsl}[3]
{\fftfuncl{newfftS}{#1}{#2}{#3}}

\newcommand{\rnsplitfftdsl}[3]{\fftfuncl{rnewfftS}{#1}{#2}{#3}}

\newcommand{\rnsplitfftdslT}[3]{\fftfunclT{rnewfftS}{#1}{#2}{#3}}

\newcommand{\splitdctII}[2]{\fftfunc{splitdctII}{#1}{#2}}

\newcommand{\nsplitdctII}[2]{\fftfunc{newdctII}{#1}{#2}}
\newcommand{\nsplitdctIII}[2]{\fftfunc{newdctIII}{#1}{#2}}

\newcommand{\indegree}{\operatorname{\mathrm{indegree}}}

\usepackage{babel}
\makeatother
\begin{document}

\title{Type-II/III DCT/DST algorithms \\
with reduced number of arithmetic operations}

\author{Xuancheng Shao and Steven G. Johnson{*}
\thanks{* Department of Mathematics, Massachusetts Institute of Technology,
Cambridge MA 02139.}
}

\maketitle

\begin{abstract}
We present algorithms for the discrete cosine transform (DCT) and
discrete sine transform (DST), of types II and III, that achieve a
lower count of real multiplications and additions than previously
published algorithms, without sacrificing numerical accuracy. Asymptotically,
the operation count is reduced from $2N\log_{2}N+O(N)$ to $\frac{17}{9}N\log_{2}N+O(N)$
for a power-of-two transform size $N$. Furthermore, we show that
an additional $N$ multiplications may be saved by a certain rescaling
of the inputs or outputs, generalizing a well-known technique for
$N=8$ by Arai~et al. These results are derived by considering the
DCT to be a special case of a DFT of length $4N$, with certain symmetries,
and then pruning redundant operations from a recent improved fast
Fourier transform algorithm (based on a recursive rescaling of the
conjugate-pair split radix algorithm). The improved algorithms for
the DCT-III, DST-II, and DST-III follow immediately from the improved
count for the DCT-II.
\end{abstract}


\section{Introduction}

In this paper, we describe recursive algorithms for the type-II and
type-III discrete cosine and sine transforms (DCT-II and DCT-III,
and also DST-II and DST-III), of power-of-two sizes, that require
fewer total real additions and multiplications than previously published
algorithms (with an asymptotic reduction of about 5.6\%), without
sacrificing numerical accuracy. Our DCT and DST algorithms are based
on a recently published fast Fourier transform (FFT) algorithm, which
reduced the operation count for the discrete Fourier transform (DFT)
compared to the previous-best split-radix algorithm \cite{Johnson07}.
The gains in this new FFT algorithm, and consequently in the new DCTs
and DSTs, stem from a recursive rescaling of the internal multiplicative
factors within an algorithm called a {}``conjugate-pair'' split-radix
FFT \cite{Kamar89,Gopinath89,Qian90,Krot92} so as to simplify some
of the multiplications. In order to derive a DCT algorithm from this
FFT, we simply consider the DCT-II to be a special case of a DFT with
real input of a certain even symmetry, and discard the redundant operations
\cite{Vetterli84,Duhamel86,DuhamelVe90,VuducDe00,FFTW05,Johnson07}.
The DCT-III, DST-II, and DST-III have identical operation counts to
the DCT-II of the same size, since the algorithms are related by simple
transpositions, permutations, and sign flips \cite{Wang82,Chan90,Lee94}.

Since 1968, the lowest total count of real additions and multiplications,
herein called {}``flops'' (floating-point operations), for the DFT
of a power-of-two size $N=2^{m}$ was achieved by the split-radix
algorithm, with $4N\log_{2}N-6N+8$ flops for $N>1$ \cite{Yavne68,Duhamel84,Martens84,Vetterli84,DuhamelVe90}.
This count was recently surpassed by new algorithms achieving a flop
count of $\frac{34}{9}N\log_{2}N+O(N)$ \cite{Johnson07,Lundy07}.
Similarly, the lowest-known flop count for the DCT-II of size $N=2^{m}>1$
was previously $2N\log_{2}N-N+2$ for a unitary normalization (with
the additive constant depending on normalization) \cite{Lee84,Vetterli84,Wang85,Duhamel86,Suehiro86,Hou87,Chan90,Arguello90,Lee94,Kok97,Takala00,Puschel03,Puschel03:dct},
and could be achieved by starting with the split-radix FFT and discarding
redundant operations \cite{Vetterli84,Duhamel86}. (Many DCT algorithms
with an unreported or larger flop count have also been described \cite{Ahmed74,Haralick76,Chen77,Narasimha78,Tseng78,Yaroslavskii79,Makhoul80,Malvar87,Li91,Steidl91,Feig92,Astola99,Guo01,Plonka05}.)
Based on our new FFT algorithm, the flop counts for the various DCT
types were reduced using a code generator \cite{Frigo99,FFTW05} that
automatically pruned the redundant operations from an FFT with a given
symmetry, but neither an explicit algorithm nor a general formula
for the flop count were presented except for DCT-I \cite{Johnson07}.
In this paper, we use the same starting point to {}``manually''
derive a DCT-II algorithm by pruning redundant operations from a real-even
FFT, and give the general formula for the new flop count (for $N=2^{m}>1$):\begin{multline}
\frac{17}{9}N\log_{2}N-\frac{17}{27}N-\frac{1}{9}(-1)^{\log_{2}N}\log_{2}N\\
+\frac{7}{54}(-1)^{\log_{2}N}+\frac{3}{2}.\label{eq:newDCT-flops}\end{multline}
The first savings over the previous record occur for $N=16$, and
are summarized in Table~\ref{tab:counts} for several $N$.%
\begin{table}[t]
\begin{centering}\begin{tabular}{|r|r|r|}
\hline 
$N$&
previous best DCT-II&
New algorithm\tabularnewline
\hline
\hline 
16&
114&
112\tabularnewline
\hline 
32&
290&
284\tabularnewline
\hline 
64&
706&
686\tabularnewline
\hline 
128&
1666&
1614\tabularnewline
\hline 
256&
3842&
3708\tabularnewline
\hline 
512&
8706&
8384\tabularnewline
\hline 
1024&
19458&
18698\tabularnewline
\hline 
2048&
43010&
41266\tabularnewline
\hline 
4096&
94210&
90264\tabularnewline
\hline
\end{tabular}\par\end{centering}

\caption{\label{tab:counts}Flop counts (real adds + mults) of previous best
DCT-II and our new algorithm }
\end{table}
 We also consider the effect of normalization on this flop count:
the above count was for a unitary transform, but slightly different
counts are obtained by other choices. Moreover, we show that a further
$N$ multiplications can be saved by individually rescaling \emph{every}
output of the DCT-II (or input of the DCT-III). In doing so, we generalize
a result by Arai \emph{et al.}, who showed that eight multiplications
could be saved by scaling the outputs of a DCT-II of size $N=8$ \cite{Arai88},
a result commonly applied to JPEG compression \cite{Pennebaker93}.

If we merely wished to show that the DCT-II/III could be computed
in $\frac{17}{9}N\log_{2}N+O(N)$ operations, we could do so by using
well-known techniques to re-express the DCT in terms of a real-input
DFT of length $N$ plus $O(N)$ pre/post-processing operations \cite{Chen77,Narasimha78,Makhoul80,Swarztrauber82,PressFlaTeu92,VanLoan92},
and then substituting our new FFT that requires $\frac{17}{9}N\log_{2}N+O(N)$
operations for real inputs \cite{Johnson07}. However, with FFT and
DCT algorithms, there is great interest in obtaining not only the
best possible asymptotic constant factor, but also the best possible
\emph{exact} count of arithmetic operations (which, for the DCT-II
of power-of-two sizes, had not changed by even one operation for over
20 years). Our result (\ref{eq:newDCT-flops}) is intended as a new
upper bound on this (still unknown) minimum exact count, and therefore
we have done our best with the $O(N)$ terms as well as the asymptotic
constant. It turns out, in fact, that our algorithm to achieve (\ref{eq:newDCT-flops})
is closely related to well-known algorithms for expressing the DCT-II
in terms of a real-input FFT of the same length, but it does not seem
obvious \emph{a priori} that this is what one obtains by pruning our
FFT for symmetric data.

In the following sections, we first review how a DCT-II may be expressed
as a special case of a DFT, and how the previous optimum flop count
can be achieved by pruning redundant operations from a conjugate-pair
split-radix FFT. Then, we briefly review the new FFT algorithm presented
in \cite{Johnson07}, and derive the new DCT-II algorithm. We follow
by considering the effect of normalization and scaling. Finally, we
consider the extension of this algorithm to algorithms for the DCT-III,
DST-II, and DST-III. We close with some concluding remarks about future
directions. We emphasize that no \emph{proven} tight lower bound on
the DCT-II flop count is currently known, and we make no claim that
eq\@.~(\ref{eq:newDCT-flops}) is the lowest possible (although
we have endeavored not to miss any obvious optimizations).

\section{DCT-II in terms of DFT\label{sec:DCT-II from DFT}}

Various forms of discrete cosine transform have been defined, corresponding
to different boundary conditions on the transform \cite{RaoYip90}.
Perhaps the most common form is the type-II DCT, used in image compression
\cite{Pennebaker93} and many other applications. The DCT-II is typically
defined as a real, orthogonal (unitary), linear transformation by
the formula (for $k=0,\ldots,N-1$):

\begin{equation}
C_{k}^{\mathrm{II}}=\sqrt{\frac{2-\delta_{k,0}}{N}}\sum_{n=0}^{N-1}x_{n}\cos\left[\frac{\pi}{N}\left(n+\frac{1}{2}\right)k\right],\label{eq:DCTII-orthogonal}\end{equation}
for $N$ inputs $x_{n}$ and $N$ outputs $C_{k}^{\mathrm{II}}$,
where $\delta_{k,0}$ is the Kronecker delta ($=1$ for $k=0$ and
$=0$ otherwise). However, we wish to emphasize in this paper that
the DCT-II (and, indeed, all types of DCT) can be viewed as special
cases of the discrete Fourier transform (DFT) with real inputs of
a certain symmetry, and where only a subset of the outputs need be
computed. This (well known) viewpoint is fruitful because it means
that any FFT algorithm for the DFT leads immediately to a corresponding
fast algorithm for the DCT-II simply by discarding the redundant operations
\cite{Vetterli84,Duhamel86,DuhamelVe90,VuducDe00,FFTW05,Johnson07}.

The discrete Fourier transform of size $N$ is defined by \begin{equation}
X_{k}=\sum_{n=0}^{N-1}x_{n}\omega_{N}^{nk},\label{eq:DFT}\end{equation}
where $\omega_{N}=e^{-\frac{2\pi i}{N}}$ is an $N$th primitive root
of unity. In order to relate this to the DCT-II, it is convenient
to choose a different normalization for the latter transform: \begin{equation}
C_{k}=2\sum_{n=0}^{N-1}x_{n}\cos\left[\frac{\pi}{N}\left(n+\frac{1}{2}\right)k\right].\label{eq:DCTII}\end{equation}
This normalization is not unitary, but it is more directly related
to the DFT and therefore more convenient for the development of algorithms.
Of course, any fast algorithm for $C_{k}$ trivially yields a fast
algorithm for $C_{k}^{\mathrm{II}}$, although the exact count of
required multiplications depends on the normalization, an issue we
discuss in more detail in section~\ref{sec:normalization}.

In order to derive $C_{k}$ from the DFT formula, one can use the
identity $2\cos(\pi\ell/N)=\omega_{4N}^{2\ell}+\omega_{4N}^{4N-2\ell}$
to write: \begin{equation}
\begin{split}C_{k} & =2\sum_{n=0}^{N-1}x_{n}\cos\left[\frac{\pi}{N}\left(n+\frac{1}{2}\right)k\right]\\
 & =\sum_{n=0}^{N-1}x_{n}\omega_{4N}^{(2n+1)k}+\sum_{n=0}^{N-1}x_{n}\omega_{4N}^{(4N-2n-1)k}\\
 & =\sum_{n=0}^{4N-1}\tilde{x}_{n}\omega_{4N}^{nk},\end{split}
\label{eq:DCTII-from-DFT}\end{equation}
where $\tilde{x}_{n}$ is a real-even sequence of length $4N$ (i.e.
$\tilde{x}_{4N-n}=\tilde{x}_{n}$), defined as follows for $0\leq n<N$:
\begin{equation}
\tilde{x}_{2n}=\tilde{x}_{4N-2n-2}=0,\label{eq:DCTII-zeros}\end{equation}
\begin{equation}
\tilde{x}_{2n+1}=\tilde{x}_{4N-(2n+1)}=x_{n}.\label{eq:DCTII-symmetry}\end{equation}
Thus, the DCT-II of size $N$ is precisely the first $N$ outputs
of a DFT of size $4N$, of real-even inputs, where the even-indexed
inputs are zero.%
\begin{figure}[t]
\begin{centering}\includegraphics[width=1\columnwidth,keepaspectratio]{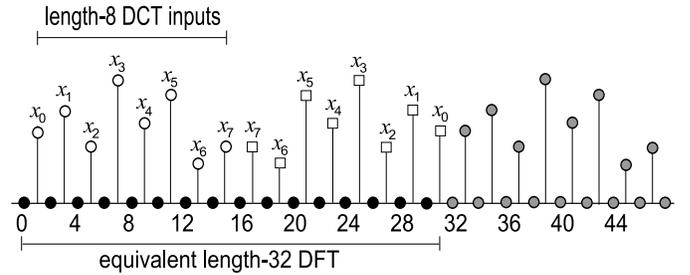}\par\end{centering}

\caption{\label{fig:dct-from-dft}A DCT-II of size $N=8$ (open dots, $x_{0},\ldots,x_{7}$)
is equivalent to a size-$4N$ DFT via interleaving with zeros (black
dots) and extending in an even (squares) periodic (gray) sequence.}
\end{figure}

This is illustrated by an example, for $N=8$, in figure~\ref{fig:dct-from-dft}.
The eight inputs of the DCT are shown as open dots, which are interleaved
with zeros (black dots) and extended to an even (squares) periodic
(gray dots) sequence of length $4N=32$ corresponding to the DFT.
(The type-II DCT is distinguished from the other DCT types by the
fact that it is even about both the left and right boundaries of the
original data, and the symmetry points fall halfway in between pairs
of the original data points.) We will refer, below, to this figure
in order to illustrate what happens when an FFT algorithm is applied
to this real-symmetric zero-interleaved data.

\section{\label{sec:conjugate-pair}Conjugate-pair FFT and DCT-II}

Although the previous minimum flop count for DCT-II algorithms can
be derived from the ordinary split-radix FFT algorithm \cite{Vetterli84,Duhamel86}
(and it can also be derived in other ways \cite{Lee84,Wang85,Suehiro86,Hou87,Chan90,Arguello90,Lee94,Kok97,Takala00,Puschel03,Puschel03:dct}),
here we will do the same thing using a variant dubbed the {}``conjugate-pair''
split-radix FFT. This algorithm was originally proposed to reduce
the number of flops \cite{Kamar89}, but was later shown to have the
same flop count as ordinary split-radix \cite{Gopinath89,Qian90,Krot92}.
It turns out, however, that the conjugate-pair algorithm exposes symmetries
in the multiplicative factors that can be exploited to reduce the
flop count by an appropriate rescaling \cite{Johnson07}, which we
will employ in the following sections.

\subsection{Conjugate-pair split-radix FFT}

Starting with the DFT of equation~(\ref{eq:DFT}), the decimation-in-time
conjugate-pair FFT splits it into three smaller DFTs: one of size
$N/2$ of the even-indexed inputs, and two of size $N/4$:\begin{multline}
X_{k}=\sum_{n_{2}=0}^{N/2-1}\omega_{N/2}^{n_{2}k}x_{2n_{2}}+\omega_{N}^{k}\sum_{n_{4}=0}^{N/4-1}\omega_{N/4}^{n_{4}k}x_{4n_{4}+1}\\
+\omega_{N}^{-k}\sum_{n_{4}=0}^{N/4-1}\omega_{N/4}^{n_{4}k}x_{4n_{4}-1},\label{eq:conjugate-pair}\end{multline}
where the indices $4n_{4}\pm1$ are computed modulo $N$. {[}In contrast,
the ordinary split-radix FFT uses $x_{4n_{4}+3}$ for the third sum
(a cyclic shift of $x_{4n_{4}-1}$), with a corresponding multiplicative
{}``twiddle'' factor of $\omega_{N}^{3k}$.] This decomposition
is repeated recursively (until base cases of size $N=1$ or $N=2$,
not shown, are reached), as shown by the pseudo-code in Algorithm~\ref{alg:simple-conjugate-pair}.
\begin{algorithm}[t]
\begin{algorithmic}

\item[\textbf{function}] $X_{k=0..N-1} \leftarrow \splitfft{N}{x_n}$:

\STATE $U_{k_2=0 \ldots N/2-1} \leftarrow \splitfft{N/2}{x_{2n_2}}$

\STATE $Z_{k_4=0 \ldots N/4-1} \leftarrow \splitfft{N/4}{x_{4n_4+1}}$

\STATE $Z'_{k_4=0 \ldots N/4-1} \leftarrow \splitfft{N/4}{x_{4n_4-1}}$

\FOR{$k=0$ to $N/4-1$}

\STATE $X_k \leftarrow U_k + \left( \omega_N^k Z_k + \omega_N^{-k} Z'_k \right)$

\STATE $X_{k+N/2} \leftarrow U_k - \left(\omega_N^k Z_k + \omega_N^{-k} Z'_k \right)$

\STATE $X_{k+N/4} \leftarrow U_{k+N/4} - i \left(\omega_N^k Z_k - \omega_N^{-k} Z'_k \right)$

\STATE $X_{k+3N/4} \leftarrow U_{k+N/4} + i \left(\omega_N^k Z_k - \omega_N^{-k} Z'_k \right)$

\ENDFOR

\end{algorithmic}

\caption{\label{alg:simple-conjugate-pair}Standard conjugate-pair split-radix
FFT algorithm of size $N$. Special-case optimizations for $k=0$
and $k=N/8$, as well as the base cases, are omitted for simplicity.}
\end{algorithm}
Here, we denote the results of the three sub-transforms of size $N/2$,
$N/4$, and $N/4$ by $U_{k}$, $Z_{k}$, and $Z_{k}'$, respectively.
The number of flops required by this algorithm, after certain simplifications
(common subexpression elimination and constant folding, and special-case
simplifications of the constants for $k=0$ and $k=N/8$) and not
counting data-independent operations like the computation of $\omega_{N}^{k}$,
is $4N\log_{2}N-6N+8$, identical to ordinary split radix. 

In the following sections, we will have to exploit further simplifications
for the case where the inputs $x_{n}$ are real. In this case, $X_{k}=X_{-k}^{*}$
(where $\ast$ denotes complex conjugation) and one can save slightly
more than half of the flops, both in the ordinary split-radix \cite{Duhamel86,SorensenJo87}
and in the conjugate-pair split-radix \cite{Johnson07}, by eliminating
redundant operations, to achieve a flop count of $2N\log_{2}N-4N+6$.
Specifically, one need only compute the outputs for $0\leq k\leq N/2$
(where the $k=0$ and $k=N/2$ outputs are purely real), and the corresponding
algorithm is shown in Algorithm~\ref{alg:rconjugate-pair}.%
\begin{algorithm}[t]
\begin{algorithmic}

\item[\textbf{function}] $X_{k=0..N/2} \leftarrow \rsplitfft{N}{x_n}$:

\STATE $U_{k_2=0 \ldots N/4} \leftarrow \rsplitfft{N/2}{x_{2n_2}}$

\STATE $Z_{k_4=0 \ldots N/8} \leftarrow \rsplitfft{N/4}{x_{4n_4+1}}$

\STATE $Z'_{k_4=0 \ldots N/8} \leftarrow \rsplitfft{N/4}{x_{4n_4-1}}$

\FOR{$k=0$ to $N/8$}

\STATE $X_k \leftarrow U_k + \left( \omega_N^k Z_k + \omega_N^{-k} Z'_k \right)$

\STATE $X_{k+N/4} \leftarrow U_{N/4-k}^{\ast} - i \left(\omega_N^k Z_k - \omega_N^{-k} Z'_k \right)$

\STATE $X_{N/4-k} \leftarrow U_{N/4-k} - i \left( \omega_N^{k} Z_{k} - \omega_N^{-k} Z_{k}' \right)^{\ast}$

\STATE $X_{N/2-k} \leftarrow U_{k}^{\ast} - \left(\omega_N^{k} Z_{k} + \omega_N^{-k} Z_{k}' \right)^{\ast}$

\ENDFOR

\end{algorithmic}

\caption{\label{alg:rconjugate-pair}Standard conjugate-pair split-radix FFT
algorithm of size $N$ as in Algorithm~\ref{alg:simple-conjugate-pair},
but specialized for the case of real inputs. Special-case optimizations
for $k=0$ and $k=N/8$, as well as the base cases, are omitted for
simplicity.}
\end{algorithm}
 Again, for simplicity we do not show special-case optimizations for
$k=0$ and $k=N/8$, so it may appear that the $X_{N/8}$, $X_{N/4}$,
and $X_{3N/8}$ outputs are computed twice. To derive Algorithm~\ref{alg:rconjugate-pair},
one exploits two facts. First, the sub-transforms operate on real
data and therefore have conjugate-symmetric output. Second, the twiddle
factors in Algorithm~\ref{alg:simple-conjugate-pair} have some redundancy
because of the identity $\omega_{N}^{\pm(N/4-k)}=\mp i\omega_{N}^{\mp k}$,
which allows us to share twiddle factors between $k$ and $N/4-k$
in the original loop.

\subsection{Fast DCT-II via split-radix FFT\label{sec:Fast-DCT-II}}

Before we derive our fast DCT-II with a reduced flop count, we first
derive a DCT-II algorithm with the \emph{same} flop count as previously
published algorithms, but starting with the conjugate-pair split-radix
algorithm. This algorithm will then be modified in section~\ref{sec:newDCT2},
below, to reduce the number of multiplications.

\begin{figure}[t]
\begin{centering}\includegraphics[width=0.8\columnwidth]{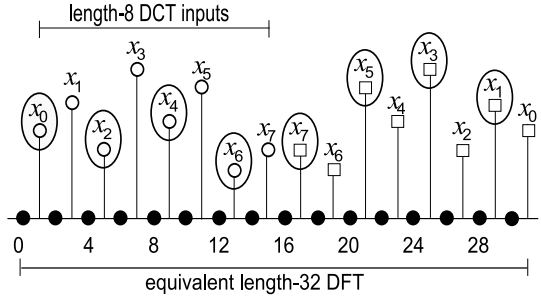}\par\end{centering}
\caption{\label{fig:splitradix-dctII}The DCT-II of $N=8$ points $x_{n}$
(open dots) is computed, in a conjugate-pair FFT of the $\tilde{N}=32$
extended data $\tilde{x}_{n}$ (squares) from figure~\ref{fig:dct-from-dft},
via the DFT $Z_{k}$ of the circled points $\tilde{x}_{4n+1}$, corresponding
to the even-indexed $x_{n}$ followed by the odd-indexed $x_{n}$.}
\end{figure}

In eq\@.~(\ref{eq:DCTII-from-DFT}), we showed that a DCT-II of
size $N$, with inputs $x_{n}$ and outputs $C_{k}$, can be expressed
as the first $N$ outputs of a DFT of size $\tilde{N}=4N$ of real-even
inputs $\tilde{x}_{n}$. {[}Therefore, $N$ in eq.~(\ref{eq:conjugate-pair})
above and in the surrounding text are replaced below by $\tilde{N}$.]
Now, consider what happens when we evaluate this DFT via the conjugate-pair
split decomposition in eq.~(\ref{eq:conjugate-pair}) and Algorithm~\ref{alg:rconjugate-pair}.
In this algorithm, we compute three smaller DFTs. First, $U_{k}$
is the DFT of size $\tilde{N}/2=2N$ of the $\tilde{x}_{2n}$ inputs,
but these are all zero and so $U_{k}$ is \emph{zero}. Second, $Z{}_{k}$
is the DFT of size $\tilde{N}/4=N$ of the $\tilde{x}_{4n+1}$ inputs
(where $4n+1$ is evaluated modulo $\tilde{N}$), which by eq\@.~(\ref{eq:DCTII-symmetry})
correspond to the original data $x_{n}$ by the formula: \begin{equation}
\tilde{y}_{n}=\tilde{x}_{4n+1}=\left\{ \begin{array}{cc}
x_{2n} & 0\leq n<N/2\\
x_{2N-1-2n} & N/2\leq n<N\end{array}\right.,\label{eq:y-tilde}\end{equation}
 where we have denoted them by $\tilde{y}_{n}$ ($0\leq n<N$) for
convenience. That is, $Z{}_{k}$ is the \emph{real-input} DFT of the
even-indexed elements of $x_{n}$ followed by the odd-indexed elements
of $x_{n}$ in reverse order. For example, this is shown for $N=8$
in figure~\ref{fig:splitradix-dctII} with the circled points corresponding
to the $\tilde{y}_{n}$, which are clearly the even-indexed $x_{n}$
followed by the odd-indexed $x_{n}$ in reverse.%
 The second DFT, $Z_{k}'$, of size $\tilde{N}/4=N$ is redundant:
it is the DFT of $\tilde{x}_{4n-1}$, but by the even symmetry of
$\tilde{x}_{n}$ this is equal to $\tilde{x}_{4(-n)+1}$, and therefore
$Z_{k}'=Z_{k}^{*}$ (the complex conjugate of $Z_{k}$).

Therefore, a fast DCT-II is obtained merely by computing a \emph{single}
real-input DFT of size $N$ to obtain $Z_{k}=Z_{N-k}^{*}$, and then
combining this according to Algorithm~\ref{alg:rconjugate-pair}
to obtain $C_{k}=X_{k}$ for $k=0\ldots N-1$. In particular, in the
loop of Algorithm~\ref{alg:rconjugate-pair}, all but the $X_{k}$
and $X_{\tilde{N}/4-k}$ terms correspond to subscripts $\geq N$
and are not needed. Moreover, we obtain \begin{align*}
X_{k} & =\omega_{\tilde{N}}^{k}Z_{k}+\omega_{\tilde{N}}^{-k}Z_{k}^{*}=2\Re\left(\omega_{4N}^{k}Z_{k}\right),\\
X_{N-k} & =-i\left(\omega_{\tilde{N}}^{k}Z_{k}-\omega_{\tilde{N}}^{-k}Z_{k}^{*}\right)^{*}=-2\Im\left(\omega_{4N}^{k}Z_{k}\right).\end{align*}
The resulting algorithm, including the special-case optimizations
for $k=0$ (where $\omega_{4N}^{0}=1$ and $Z_{0}$ is real) and $k=N/2$
(where $\omega_{4N}^{N/2}=(1-i)/\sqrt{2}$ and $Z_{N/2}$ is real),
is shown in Algorithm~\ref{alg:splitdctII}.%
\begin{algorithm}[t]
\begin{algorithmic}

\item[\textbf{function}] $C_{k=0..N-1} \leftarrow \splitdctII{N}{x_n}$:

\FOR{$n=0$ to $N/2-1$}

\STATE $\tilde{y}_n \leftarrow x_{2n}$

\STATE $\tilde{y}_{N-1-n} \leftarrow x_{2n+1}$

\ENDFOR

\STATE $Z_{k=0 \ldots N/2} \leftarrow \rsplitfft{N}{\tilde{y}_n}$

\STATE $C_0 \leftarrow 2 Z_0$

\FOR{$k=1$ to $N/2-1$}

\STATE $C_k \leftarrow 2\Re \left( \omega_{4N}^k Z_k \right)$

\STATE $C_{N-k} \leftarrow -2 \Im \left(\omega_{4N}^k Z_k \right)$

\ENDFOR

\STATE $C_{N/2} \leftarrow \sqrt{2} Z_{N/2}$

\end{algorithmic}

\caption{\label{alg:splitdctII}Fast DCT-II algorithm, matching previous best
flop count, derived from Algorithm~\ref{alg:rconjugate-pair} by
discarding redundant operations.}
\end{algorithm}

In fact, Algorithm~\ref{alg:splitdctII} is equivalent to an algorithm
derived in a different way by various authors~\cite{Narasimha78,Makhoul80}
to express a DCT-II in terms of a real-input DFT of the same length.
Here, we see that this algorithm is exactly equivalent to a conjugate-pair
split-radix FFT of length $4N$. Previous authors obtained a suboptimal
flop count $\frac{5}{2}N\log_{2}N+O(N)$ for the DCT-II using this
algorithm \cite{Makhoul80} only because they used a suboptimal real-input
DFT for the subtransform (split-radix being then almost unknown).
Using a real-input split-radix DFT to compute $Z_{k}$, the flop count
for $Z_{k}$ is $2N\log_{2}N-4N+6$ from above. To get the total flop
count for Algorithm~\ref{alg:splitdctII}, we need to add $N/2-1$
general complex multiplications by $2\omega_{4N}^{k}$ (6 flops each)
plus two real multiplications (2 flops), for a total of $2N\log_{2}N-N+2$
flops. This matches the best-known flop count in the literature (where
the $+2$ can be removed merely by choosing a different normalization
as discussed in section~\ref{sec:normalization}).

\section{New FFT and DCT-II}

In this section, we first review the new FFT algorithm introduced
in our previous work based on a recursive rescaling of the conjugate-pair
FFT \cite{Johnson07}, and then apply it to derive a fast DCT-II algorithm
as in the previous section.

\subsection{New FFT\label{sec:New-FFT}}

Based on the conjugate-pair split-radix FFT from section~\ref{sec:conjugate-pair},
a new FFT algorithm with a reduced number of flops can be derived
by scaling the subtransforms \cite{Johnson07}. We will not reproduce
the derivation here, but will simply summarize the results. In particular,
the original conjugate-pair split-radix Algorithm~\ref{alg:simple-conjugate-pair}
is split into four mutually recursive subroutines, each of which has
the same split-radix structure but computes a DFT scaled by a different
factor. These algorithms are shown in Algorithm~\ref{cap:nsplitfft}
for the case of complex inputs, and in Algorithm~\ref{cap:rnsplitfft}
specialized for real inputs, in both of which the scaling factors
are combined with the twiddle factors $\omega_{N}^{k}$ to reduce
the total number of multiplications compared to Algorithms~\ref{alg:simple-conjugate-pair}
and~\ref{alg:rconjugate-pair}, respectively.%
\begin{algorithm}[t]
\begin{algorithmic}

\item[\textbf{function}] $X_{k=0..N-1} \leftarrow \nsplitfftdsl{N}{\ell}{x_n}$:

\COMMENT{computes DFT / $s_{\ell N,k}$ for $\ell = 0,1,2$}

\STATE $U_{k_2=0 \ldots N/2-1} \leftarrow \nsplitfftdsl{N/2}{2\ell}{x_{2n_2}}$

\STATE $Z_{k_4=0 \ldots N/4-1} \leftarrow \nsplitfftdsl{N/4}{1}{x_{4n_4+1}}$

\STATE $Z'_{k_4=0 \ldots N/4-1} \leftarrow \nsplitfftdsl{N/4}{1}{x_{4n_4-1}}$

\FOR{$k=0$ to $N/4-1$}

\STATE $X_k \leftarrow U_k + \left( t_{N,k} Z_k + t_{N,k}^* Z'_k \right) \cdot (s_{N,k} / s_{\ell N,k})$

\STATE $X_{k+N/2} \leftarrow U_k - \left(t_{N,k} Z_k + t_{N,k}^* Z'_k \right)\cdot (s_{N,k} / s_{\ell N,k})$

\STATE $X_{k+N/4} \leftarrow U_{k+N/4} $ \\ \hspace{1em} $ - \: i \left(t_{N,k} Z_k - t_{N,k}^* Z'_k \right)\cdot (s_{N,k} / s_{\ell N,k+N/4})$

\STATE $X_{k+3N/4} \leftarrow U_{k+N/4} $ \\ \hspace{1em} $ + \: i \left(t_{N,k} Z_k - t_{N,k}^* Z'_k \right)\cdot (s_{N,k} / s_{\ell N,k+N/4})$

\ENDFOR

\vskip 3pt

\item[\textbf{function}] $X_{k=0..N-1} \leftarrow \nsplitfftdsfour{N}{x_n}$:

\COMMENT{computes DFT / $s_{4N,k}$}

\STATE $U_{k_2=0 \ldots N/2-1} \leftarrow \nsplitfftdsl{N/2}{2}{x_{2n_2}}$

\STATE $Z_{k_4=0 \ldots N/4-1} \leftarrow \nsplitfftdsl{N/4}{1}{x_{4n_4+1}}$

\STATE $Z'_{k_4=0 \ldots N/4-1} \leftarrow \nsplitfftdsl{N/4}{1}{x_{4n_4-1}}$

\FOR{$k=0$ to $N/4-1$}

\STATE $X_k \leftarrow \left[ U_k + \left( t_{N,k} Z_k + t_{N,k}^* Z'_k \right) \right] \cdot (s_{N,k} / s_{4N,k})$

\STATE $X_{k+N/2} \leftarrow \left[ U_k - \left(t_{N,k} Z_k + t_{N,k}^* Z'_k \right) \right] $ \\ \hspace{5em} $ \cdot \: (s_{N,k} / s_{4N,k+N/2})$

\STATE $X_{k+N/4} \leftarrow \left[ U_{k+N/4} - i \left(t_{N,k} Z_k - t_{N,k}^* Z'_k \right) \right] $ \\ \hspace{5em} $ \cdot \: (s_{N,k} / s_{4N,k+N/4})$

\STATE $X_{k+3N/4} \leftarrow \left[ U_{k+N/4} + i \left(t_{N,k} Z_k - t_{N,k}^* Z'_k \right) \right] $ \\ \hspace{5em} $ \cdot \: (s_{N,k} / s_{4N,k+3N/4})$

\ENDFOR

\end{algorithmic}

\caption{\label{cap:nsplitfft}New FFT algorithm \cite{Johnson07} of length
$N$ (divisible by $4$). The sub-transforms $\nsplitfftdsl N\ell x$
for $\ell\neq0$ are scaled by $s_{\ell N,k}$, respectively, while
$\ell=0$ is the final unscaled DFT ($s_{0,k}=1$). Special-case optimizations
for $k=0$ and $k=N/8$ are not shown.}
\end{algorithm}
\begin{algorithm}[t]
\begin{algorithmic}

\item[\textbf{function}] $X_{k=0..N/2} \leftarrow \rnsplitfftdsl{N}{\ell}{x_n}$:

\COMMENT{computes DFT / $s_{\ell N,k}$ for $\ell = 0,1,2$}

\STATE $U_{k_2=0 \ldots N/4} \leftarrow \rnsplitfftdsl{N/2}{2\ell}{x_{2n_2}}$

\STATE $Z_{k_4=0 \ldots N/8} \leftarrow \rnsplitfftdsl{N/4}{1}{x_{4n_4+1}}$

\STATE $Z'_{k_4=0 \ldots N/8} \leftarrow \rnsplitfftdsl{N/4}{1}{x_{4n_4-1}}$

\FOR{$k=0$ to $N/8$}

\STATE $X_k \leftarrow U_k + \left( t_{N,k} Z_k + t_{N,k}^* Z'_k \right) \cdot (s_{N,k} / s_{\ell N,k})$

\STATE $X_{k+N/4} \leftarrow U_{N/4-k}^{\ast} $ \\ \hspace{1em} $ - \: i \left(t_{N,k} Z_k - t_{N,k}^* Z'_k \right)\cdot (s_{N,k} / s_{\ell N,k+N/4})$

\STATE $X_{N/4-k} \leftarrow U_{N/4-k} $ \\ \hspace{1em} $ - \: i \left( t_{N,k} Z_k - t_{N,k}^* Z_k' \right)^* \cdot (s_{N,k} / s_{\ell N,k+N/4})$

\STATE $X_{N/2-k} \leftarrow U_k^* $ \\ \hspace{1em} $ - \left(t_{N,k} Z_k + t_{N,k}^* Z_k' \right)^* \cdot (s_{N,k} / s_{\ell N, k})$

\ENDFOR

\vskip 3pt

\item[\textbf{function}] $X_{k=0..N/2} \leftarrow \rnsplitfftdsfour{N}{x_n}$:

\COMMENT{computes DFT / $s_{4N,k}$}

\STATE $U_{k_2=0 \ldots N/4} \leftarrow \rnsplitfftdsl{N/2}{2}{x_{2n_2}}$

\STATE $Z_{k_4=0 \ldots N/8} \leftarrow \rnsplitfftdsl{N/4}{1}{x_{4n_4+1}}$

\STATE $Z'_{k_4=0 \ldots N/8} \leftarrow \rnsplitfftdsl{N/4}{1}{x_{4n_4-1}}$

\FOR{$k=0$ to $N/8$}

\STATE $X_k \leftarrow \left[ U_k + \left( t_{N,k} Z_k + t_{N,k}^* Z'_k \right) \right] \cdot (s_{N,k} / s_{4N,k})$

\STATE $X_{k+N/4} \leftarrow \left[ U_{N/4-k}^* - i \left(t_{N,k} Z_k - t_{N,k}^* Z'_k \right) \right] $ \\ \hspace{5em} $ \cdot \: (s_{N,k} / s_{4N,k+N/4})$

\STATE $X_{N/4-k} \leftarrow \left[ U_{N/4-k} - i \left(t_{N,k} Z_k - t_{N,k}^* Z_k' \right)^* \right] $ \\ \hspace{5em} $ \cdot \: (s_{N,k} / s_{4N,N/4-k})$

\STATE $X_{N/2-k} \leftarrow \left[ U_k^* - \left(t_{N,k} Z_k + t_{N,k}^* Z_k' \right)^* \right] $ \\ \hspace{5em} $ \cdot \: (s_{N,k} / s_{4N,N/2-k})$

\ENDFOR

\end{algorithmic}

\caption{\label{cap:rnsplitfft}New FFT algorithm of size $N$ (divisible
by $4$), as in Algorithm~\ref{cap:nsplitfft} but specialized for
the case of real inputs. The sub-transforms $\rnsplitfftdsl N\ell x$
for $\ell\neq0$ are scaled by $s_{\ell N,k}$, respectively, while
$\ell=0$ is the final unscaled DFT ($s_{0,k}=1$). Special-case optimizations
for $k=0$ and $k=N/8$ are not shown.}
\end{algorithm}
 Again, special cases for $k=0$ and $k=N/8$, as well as the $N=1,2$
base cases and obvious simplifications such as common-subexpression
elimination, are not shown for simplicity.

The key aspect of these algorithms is the scale factor $s_{N,k}$,
where the subtransforms compute the DFT scaled by $1/s_{\ell N,k}$
for $\ell=1,2,4$. (The $\ell=0,1,2$ subroutines are condensed into
one in Algorithms~\ref{cap:nsplitfft}--\ref{cap:rnsplitfft}, by
making $\ell$ a variable, but in practice they need to be implemented
separately in order to exploit the special cases of the scale factor
\cite{Johnson07}. The $\ell=4$ case is written separately because
it is factorized differently.) This scale factor is defined for $N=2^{m}$
by the following recurrence, where $k_{4}=k$ mod $\frac{N}{4}$:

\begin{equation}
s_{N=2^{m},k}=\left\{ \begin{array}{cl}
1 & \quad\textrm{for }N\leq4\\
s_{N/4,k_{4}}\cos(2\pi k_{4}/N) & \quad\textrm{for }k_{4}\leq N/8\\
s_{N/4,k_{4}}\sin(2\pi k_{4}/N) & \quad\textrm{otherwise}\end{array}\right..\label{eq:scalesym}\end{equation}
 This definition has the properties: $s_{N,0}=1$, $s_{N,k+N/4}=s_{N,k}$,
and $s_{N,N/4-k}=s_{N,k}$. Also, $s_{N,k}>0$ and decays rather slowly
with $N$: $s_{N,k}$ has an $\Omega(N^{\log_{4}\cos(\pi/5)})$ asymptotic
lower bound \cite{Johnson07}. When these scale factors are combined
with the twiddle factors $\omega_{N}^{k}$, we obtain terms of the
form \begin{equation}
t_{N,k}=\omega_{N}^{k}\frac{s_{N/4,k}}{s_{N,k}},\label{eq:t}\end{equation}
which is always a complex number of the form $\pm1\pm i\tan\frac{2\pi k}{N}$
or $\pm\cot\frac{2\pi k}{N}\pm i$ and can therefore be multiplied
with two fewer real multiplications than are required to multiply
by $\omega_{N}^{k}$. Because of the symmetry $s_{N,N/4-k}=s_{N,k}$,
it is possible to specialize Algorithm~\ref{cap:nsplitfft} for real
data in the same way as in Algorithm~\ref{alg:rconjugate-pair},
because the scale factor preserves the conjugate symmetry $X_{k}=X_{N-k}^{*}$
of the outputs of all subtransforms \cite{Johnson07}.%
\footnote{In Algorithm~\ref{cap:rnsplitfft}, we have also utilized various
symmetries of the scale factors, such as $s_{2N,N/4-k}=s_{2N,k+N/4}$,
as described in our previous work \cite{Johnson07}, to make explicit
the shared multiplicative factors between the different terms.%
}

The resulting flop count, for arbitrary complex data $x_{n}$, is
then reduced from $4N\log_{2}N-6N+8$ to\begin{multline}
\frac{34}{9}N\log_{2}N-\frac{124}{27}N-2\log_{2}N\\
-\frac{2}{9}(-1)^{\log_{2}N}\log_{2}N+\frac{16}{27}(-1)^{\log_{2}N}+8.\label{eq:newsplit-flops}\end{multline}
In particular, assuming that complex multiplications are implemented
as four real multiplications and two real additions, then the savings
are purely in the number of real multiplications. The number $M(N)$
of real multiplications saved over ordinary split radix is given by:\begin{multline}
M(N)=\frac{2}{9}N\log_{2}N-\frac{38}{27}N+2\log_{2}N\\
+\frac{2}{9}(-1)^{\log_{2}N}\log_{2}N-\frac{16}{27}(-1)^{\log_{2}N}.\label{eq:M}\end{multline}
If the data are purely real, it was shown that $M(N)/2$ multiplications
are saved over the corresponding real-input split-radix algorithm
\cite{Johnson07}. These flop counts are to compute the original,
unscaled definition of the DFT. If one is allowed to scale the outputs
by any factor desired, then scaling by $1/s_{N,k}$ {[}corresponding
to $\nsplitfftdsl{N}{1}{x_{n}}$ in Algorithm~\ref{cap:nsplitfft}],
saves an additional $M_{S}(N)-M(N)\geq0$ multiplications for complex
inputs, where: \begin{multline}
M_{S}(N)=\frac{2}{9}N\log_{2}N-\frac{20}{27}N\\
+\frac{2}{9}(-1)^{\log_{2}N}\log_{2}N-\frac{7}{27}(-1)^{\log_{2}N}+1.\label{eq:Ms}\end{multline}
 For real inputs, one similarly saves $M_{S}(N)/2$ flops for $\rnsplitfftdsl{N}{1}{x_{n}}$
operating on real $x_{n}$ in Algorithm~\ref{cap:rnsplitfft}.

Although the division by a cosine or sine in the scale factor may,
at first glance, seem as if it may exacerbate floating-point errors,
this is not the case. Cooley--Tukey based FFT algorithms have $O(\sqrt{\log N})$
root-mean-square error growth, and $O(\log N)$ error bounds, in finite-precision
floating-point arithmetic (assuming that the trigonometric constants
are precomputed accurately) \cite{GenSan66,Schatzman96,Tasche00},
and the new FFT is no different \cite{Johnson07}. The reason for
this is straightforward: it never adds or subtracts scaled and unscaled
values. Instead, wherever the original FFT would compute $a+b$, the
new FFT computes $s\cdot(a+b)$ for some fixed scale factor $s$.

\subsection{Fast DCT-II from new FFT\label{sec:newDCT2}}

To derive the new DCT-II algorithm based on the new FFT of the previous
section, we follow exactly the same process as in Sec\@.~\ref{sec:Fast-DCT-II}.
We express the DCT-II of length $N$ in terms of a DFT of length $4N$,
apply the FFT algorithm~\ref{cap:rnsplitfft}, and discard the redundant
operations. As before, the $U_{k}$ subtransform is identically zero,
$Z_{k}$ is the transform of the even-indexed elements of $x_{n}$
followed by the odd-indexed elements in reverse order, and $Z_{k}'=Z_{N-k}^{*}$
is redundant.%
\begin{algorithm}[t]
\begin{algorithmic}

\item[\textbf{function}] $C_{k=0..N-1} \leftarrow \nsplitdctII{N}{x_n}$:

\FOR{$n=0$ to $N/2-1$}

\STATE $\tilde{y}_n \leftarrow x_{2n}$

\STATE $\tilde{y}_{N-1-n} \leftarrow x_{2n+1}$

\ENDFOR

\STATE $Z_{k=0..N/2} \leftarrow \rnsplitfftdsl{N}{1}{\tilde{y}_n}$

\STATE $C_0 \leftarrow 2 Z_0$

\FOR{$k=1$ to $N/2-1$}

\STATE $C_k \leftarrow 2\Re \left( \omega_{4N}^k s_{N,k} Z_k \right)$

\STATE $C_{N-k} \leftarrow - 2\Im \left( \omega_{4N}^k s_{N,k} Z_k \right)$

\ENDFOR

\STATE $C_{N/2} \leftarrow \sqrt {2} Z_{N/2}$

\end{algorithmic}

\caption{\label{alg:nsplitdctII}New DCT-II algorithm of size $N=2^{m}$,
based on discarding redundant operations from the new FFT algorithm
of size $4N$, achieving new record flop count.}
\end{algorithm}

The resulting algorithm is shown in Algorithm~\ref{alg:nsplitdctII},
and differs from the original fast DCT-II of Algorithm~\ref{alg:splitdctII}
in only two ways. First, instead of calling the ordinary split-radix
(or conjugate-pair) FFT for the subtransform, it calls $\nsplitfftds{N}{x_{n}}$.
Second, because this subtransform is scaled by $1/s_{N,k}$, the twiddle
factor $\omega_{4N}^{k}$ in Algorithm~\ref{alg:nsplitdctII} is
multiplied by $s_{N,k}$.

To derive the flop count for Algorithm~\ref{alg:nsplitdctII}, we
merely need to add the flop count for the subtransform {[}which saves
$M_{S}(N)/2=\frac{1}{9}N\log_{2}N+O(N)$ flops, from eq\@.~(\ref{eq:Ms}),
compared to ordinary real-input split radix] with the number of flops
in the loop, where the latter is exactly the same as for Algorithm~\ref{alg:splitdctII}
because the products $\omega_{4N}^{k}s_{N,k}$ can be precomputed.
Therefore, we save a total of $M_{S}(N)/2$ flops compared to the
previous best flop count of $2N\log_{2}N-N+2$, resulting in the flop
count of eq\@.~(\ref{eq:newDCT-flops}). This formula matches the
flop count that was achieved by an automatic code generator in Ref\@.~\cite{Johnson07}.

Because the new DCT algorithm is mathematically merely a special set
of inputs for the underlying FFT, it will have the same favorable
logarithmic error bounds as discussed in the previous section.

\section{Normalizations\label{sec:normalization}}

In the above definition of DCT-II, we use a scale factor {}``2''
in front of the summation in order to make the DCT-II directly equivalent
to the corresponding DFT. But in some circumstances, it is useful
to multiply by other factors, and different normalizations lead to
slightly different flop counts. For example, one could use the unitary
normalization from eq\@.~(\ref{eq:DCTII-orthogonal}), which replaces
$2$ by $\sqrt{2/N}$ or $\sqrt{1/N}$ (for $k=0$) and requires the
same number of flops. If one uses the unitary normalization multiplied
by $\sqrt{N}$, then one saves two multiplications in the $k=0$ and
$k=N/2$ terms (in this normalization, $C_{0}=Z_{0}$ and $C_{N/2}=Z_{N/2}$)
for both Algorithm~\ref{alg:splitdctII} and Algorithm~\ref{alg:nsplitdctII}.
(This leads to a commonly quoted $2N\log_{2}N-N$ formula for the
previous flop count.)

It is also common to compute a DCT-II with scaled outputs, \emph{e.g.}
for the JPEG image-compression standard where an arbitrary scaling
can be absorbed into a subsequent quantization step~\cite{Pennebaker93},
and in this case the scaling can save 8 multiplications~\cite{Arai88}
over the 42 flops required for an unitary 8-point DCT-II. Since our
$\nsplitfftdsl{N}{1}{x_{n}}$ attempts to be the optimally scaled
FFT, we should be able to derive this scaled DCT-II by using $\rnsplitfftdsl{N}{1}{x_{n}}/2$
for our length-$4N$ DFT instead of $\rnsplitfftdsl{N}{0}x$, and
again pruning the redundant operations. Doing so, we obtain an algorithm
identical to Algorithm~\ref{alg:nsplitdctII} except that $2\omega_{4N}^{k}s_{N,k}$
is replaced by $t_{4N,k}$, which requires fewer multiplications,
and also we now obtain $C_{0}=Z_{0}$ and $C_{N/2}=Z_{N/2}$. This
algorithm computes $C_{k}/2s_{4N,k}$ instead of $C_{k}$. In particular,
we save exactly $N$ multiplications over Algorithm~\ref{alg:nsplitdctII},
which matches the result by Ref\@.~\cite{Arai88} for $N=8$ but
generalizes it to all $N$. This savings of $N$ multiplications for
a scaled DCT-II also matches the flop count that was achieved by an
automatic code generator in Ref\@.~\cite{Johnson07}.

\section{Fast DCT-III, DST-II, and DST-III}

Given any algorithm for the DCT-II, one can immediately obtain algorithms
for the DCT-III, DST-II, and DST-III with exactly the same number
of arithmetic operations. In this way, any improvement in the DCT-II,
such as the one described in this paper, immediately leads to improved
algorithms for those transforms and vice versa. In this section, we
review the equivalence between those transforms.

\subsection{DCT-III}

\begin{figure}[t]
\begin{centering}\includegraphics[width=1\columnwidth]{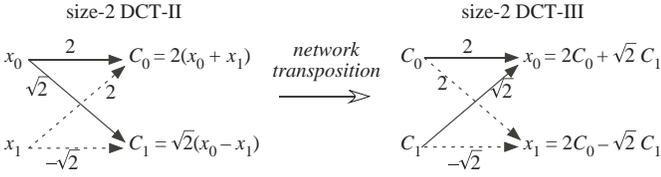}\par\end{centering}

\caption{\label{fig:network-transposition}Example of network transposition
for a size-2 DCT-II (left). When the linear network is transposed
(all edges are reversed), the resulting network computes the transposed-matrix
operation: a size-2 DCT-III (right).}
\end{figure}

To obtain a DCT-III algorithm from a DCT-II algorithm, one can exploit
two facts. First, the DCT-III, viewed as a matrix, is simply the transpose
of the DCT-II matrix. Second, any linear algorithm can be viewed as
a \emph{linear network}, and the transpose operation is computed by
the \emph{network transposition} of this algorithm~\cite{CrochiereOp75}.
To review, a linear network represents the algorithm by a directed
acyclic graph, where the edges represent multiplications by constants
and the vertices represent additions of the incoming edges. Network
transposition simply consists of reversing the direction of every
edge. We prove below that the transposed network requires the same
number of additions and multiplications for networks with the same
number of inputs and outputs, and therefore the DCT-III can be computed
in the same number of flops as the DCT-II by the transposed algorithm.
The DCT-III corresponding to the transpose of eq\@.~(\ref{eq:DCTII})
is\begin{equation}
C{}_{k}^{\mathrm{T}}=2\sum_{n=0}^{N-1}x_{n}\cos\left[\frac{\pi}{N}n\left(k+\frac{1}{2}\right)\right].\label{eq:DCTIII}\end{equation}
 Since the DCT-III is the transpose of the DCT-II, a fast DCT-III
algorithm is obtained by network transposition of a fast DCT-II algorithm.
For example, network transposition of a size-2 DCT-III is shown in
figure\@.~\ref{fig:network-transposition}.%

A proof that the transposed network requires the same number of flops
is as follows. Clearly, the number of multiplications, the number
of edges with weight $\neq\pm1$, is unchanged by transposition. The
number of additions is given by the sum of $\indegree-1$ for all
the vertices, except for the $N$ input vertices which have indegree
zero. That is, for a set $V$ of vertices and a set $E$ of edges:\begin{align}
\mbox{\# adds} & =N+\sum_{v\in V}[\indegree(v)-1]\nonumber \\
 & =N+|E|-|V|,\label{eq:graph}\end{align}
using the fact that the sum of the $\indegree$ over all vertices
is simply the number $|E|$ of edges. Because the above expression
is obviously invariant under network transposition as long as $N$
is unchanged (i.e. same number of inputs and outputs), the number
of additions is invariant.

More explicitly, a fast DCT-III algorithm derived from the transpose
of our new Algorithm~\ref{alg:nsplitdctII} is outlined in Algorithm~\ref{alg:newDCTIII}.
Whereas for the DCT-II we computed the real-input DFT (with conjugate-symmetric
output) of the rearranged inputs $\tilde{y}_{n}$ and then post-processed
the complex outputs $Z_{k}=Z_{N-k}^{*}$, now we do the reverse. We
first preprocess the inputs to form a complex array $Z_{k}=Z_{N-k}^{*}$,
then perform a real-\emph{output}, scaled-\emph{input} DFT to obtain
$\tilde{y}_{k}$, and finally assign the even and odd elements of
the result $C_{k}^{\mathrm{T}}$ from the first and second halves
of $\tilde{y}_{k}$. Without the scale factors, this is equivalent
to a well-known algorithm to express a DCT-III in terms of a realoutput
DFT of the same size \cite{Narasimha78,Makhoul80}. In order to minimize
the flop count once the scale factors are included, however, it is
important that this real-output DFT operate on scaled \emph{inputs}
rather than scaled \emph{outputs}, so it is intrinsically different
(transposed) from Algorithm~\ref{cap:nsplitfft}. The real-output
(scaled, conjugate-symmetric input) version of $\rnsplitfftdsl{N}{1}{x_{n}}$
can be derived by network transposition of the real-input case (since
network transposition interchanges inputs and outputs). Equivalently,
whereas the $\rnsplitfftdsl{N}{1}{x_{n}}$ was a scaled-output decimation-in-time
algorithm specialized for real inputs, the transposed algorithm $\rnsplitfftdslT{N}{1}{x_{n}}$
is a scaled-input decimation-in-frequency algorithm specialized for
real outputs. We do not list $\rnsplitfftdslT{N}{1}{x_{n}}$ explicitly
here, however; our main point is simply to establish that an algorithm
for the DCT-III with exactly the same flop count (\ref{eq:newDCT-flops})
follows immediately from the new DCT-II.%
\begin{algorithm}[t]
\begin{algorithmic}

\item[\textbf{function}] $D_{k=0..N-1} \leftarrow \nsplitdctIII{N}{x_n}$:

\STATE $Z_0 \leftarrow 2 x_0$

\FOR{$n=1$ to $N/2-1$}

\STATE $Z_n \leftarrow 2 \omega_{4N}^{-n} s_{N,n} (x_n - i x_{N-n})$

\STATE $Z_{N-n} \leftarrow Z_n^*$

\ENDFOR

\STATE $Z_{N/2} \leftarrow \sqrt{2} x_{N/2}$

\STATE $\tilde{y}_k \leftarrow \rnsplitfftdslT{N}{1}{Z_n}$

\FOR{$k=0$ to $N/2-1$}

\STATE $C_{2k}^\mathrm{T} \leftarrow \tilde{y}_k$

\STATE $C_{2k+1}^\mathrm{T} \leftarrow \tilde{y}_{N-1-k}$

\ENDFOR

\end{algorithmic}

\caption{\label{alg:newDCTIII}New DCT-III algorithm of size $N=2^{m}$, based
on network transposition of Algorithm~\ref{alg:nsplitdctII}, with
the same flop count.}
\end{algorithm}

There are also other ways to derive a DCT-III algorithm with the same
flop count without using network transposition, of course. For example,
one can again consider the DCT-III to be a real-input DFT of size
$4N$ with certain symmetries (different from the DCT-II's symmetries),
and prune redundant operations from Algorithm~\ref{cap:rnsplitfft},
resulting in a decimation-in-time algorithm. Such a derivation is
useful because it allows one to derive a scaled-\emph{output} DCT-III
as well, which turns out to be a subroutine of a new DCT-IV algorithm
that we describe in another manuscript, currently in press \cite{ShaoJo08-dct4}.
However, this derivation is more complicated than the one above, resulting
in four mutually recursive DCT-III subroutines, and does not lower
the flop count, so we omit it here.

\subsection{DST-II and DST-III}

The DST-II is defined, using a unitary normalization, as: \begin{equation}
S_{k}^{\mathrm{II}}=\sqrt{\frac{2-\delta_{k,N}}{N}}\sum_{n=0}^{N-1}x_{n}\sin\left[\frac{\pi}{N}\left(n+\frac{1}{2}\right)k\right]\label{eq:dstII}\end{equation}
for $k=1,\ldots,N$ (\emph{not} $0,\ldots,N-1$). As in section~\ref{sec:DCT-II from DFT},
it is more convenient to develop algorithms for an unnormalized variation:
\begin{equation}
S_{k}=2\sum_{n=0}^{N-1}x_{n}\sin\left[\frac{\pi}{N}\left(n+\frac{1}{2}\right)k\right]\label{eq:dstII-norm2}\end{equation}
similar to our normalization of $C_{k}$ and $C_{k}^{\mathrm{T}}$
above. Although we could derive fast algorithms for $S_{k}$ directly
by treating it as a DFT of length $4N$ with odd symmetry, interleaved
with zeros, and discarding redundant operations similar to above,
it turns out there is a simpler technique. The DST-II is \emph{exactly}
equivalent to a DCT-II in which the outputs are reversed and every
other input is multiplied by $-1$ \cite{Wang82,Chan90,Lee94}: \begin{equation}
S_{N-k}=2\sum_{n=0}^{N-1}(-1)^{n}x_{n}\cos\left[\frac{\pi}{N}\left(n+\frac{1}{2}\right)k\right]\label{eq:dstII-dctII}\end{equation}
for $k=0,\ldots,N-1$. It therefore follows that a DST-II can be computed
with the same number of flops as a DCT-II of the same size, assuming
that multiplications by $-1$ are free---the reason for this is that
sign flips can be absorbed at no cost by converting additions into
subtractions or vice versa in the subsequent algorithmic steps. Therefore,
our new flop count (\ref{eq:newDCT-flops}) immediately applies to
the DST-II.

Similarly, a DST-III is given by the transpose of the DST-II (which
is the inverse, for the unitary normalization). In unnormalized form,
this is \begin{equation}
S_{k}^{\mathrm{T}}=2\sum_{n=1}^{N}x_{n}\sin\left[\frac{\pi}{N}n\left(k+\frac{1}{2}\right)\right]\label{eq:dstIII}\end{equation}
for $k=0,\ldots,N-1$. This must have the same number of flops as
the DST-II by the network transposition argument above. Alternatively,
it can also be obtained from the DCT-III by reversing the inputs and
multiplying every other output by $-1$ \cite{Chan90}: \begin{equation}
S_{k}^{\mathrm{T}}=(-1)^{k}\sum_{n=0}^{N-1}x_{N-n}\cos\left[\frac{\pi}{N}n\left(k+\frac{1}{2}\right)\right],\label{eq:dstIII-dctIII}\end{equation}
which can be obtained from $C_{k}^{T}$ with the same number of flops.

\section{Conclusion}

We have shown that a improved count of real additions and multiplications
(flops), eq\@.~(\ref{eq:newDCT-flops}), can be obtained for the
DCT/DST types II/III by pruning redundant operations from a recent
improved FFT algorithm. We have also shown how to save $N$ multiplications
by rescaling the outputs of a DCT-II (or the inputs of a DCT-III),
generalizing a well-known technique for $N=8$ \cite{Arai88}. However,
we do not claim that our new count is optimal---it may be that further
gains can be obtained by, for example, extending the rescaling technique
of \cite{Johnson07} to greater generality. In particular, as has
been pointed out by other authors \cite{DuhamelVe90}, any improvement
in the arithmetic complexity or flop counts for the DFT immediately
leads to corresponding improvements in DCTs/DSTs, and vice versa.
Our most important result, we believe, is the fact that there suddenly
appears to be new room for improvement in problems that had seen no
reductions in flop counts for many years.

The question of the minimal operation counts for basic linear transformations
such as DCTs and DSTs is of longstanding theoretical interest. The
practical impact of a few percent improvement in flop counts, admittedly,
is less clear because computation time on modern hardware is often
not dominated by arithmetic \cite{FFTW05}. Nevertheless, minimal-arithmetic
hard-coded DCTs of small sizes are often used in audio and image compression
\cite{RaoYip90,Pennebaker93}, and the availability of any new algorithm
with regular structure amenable to implementation leads to rich new
opportunities for performance optimization.

Similar arithmetic improvements also occur for other types of DCT
and DST \cite{Johnson07}, as well as for the modified DCT (MDCT,
closely related to a DCT-IV), and the explicit description of those
algorithms is the subject of another manuscript currently in press
\cite{ShaoJo08-dct4}. We have already shown that the new exact count
for the DCT-I is $2N\log_{2}N-3N+2\log_{2}N+5-M(2N)/4=\frac{17}{9}N\log_{2}N+O(N)$,
where $M(N)$ is given by eq\@.~(\ref{eq:M}) \cite{Johnson07}.
However, no new exact count for arbitrary $N=2^{m}$ has yet been
published for the DCT-IV and related transforms. Again, it immediately
follows that the asymptotic cost of a DCT-IV (and hence an MDCT) is
reduced to $\frac{17}{9}N\log_{2}N+O(N)$ simply by applying known
algorithms to express a DCT-IV in terms of a DCT-II of size $N$ \cite{Chan90}
(although this algorithm has large numerical errors \cite{FFTW05})
or in terms of two DCT-II/III transforms of size $N/2$ \cite{Wang85}.
However, again we wish to establish a new lowest-known upper bound
for the \emph{exact} count of operations of the DCT-IV, and not just
the asymptotic constant factor. (It turns out that the DCT-IV algorithm
obtained by pruning redundant computations from our new FFT, this
time of symmetric data with size $8N$ \cite{FFTW05}, is closely
related to the algorithm by Wang \emph{et al}. \cite{Wang85}, just
as as the algorithm for the DCT-II in this paper was closely related
to a known algorithm derived by other means \cite{Narasimha78,Makhoul80}.)

This work was supported in part by a grant from the MIT Undergraduate
Research Opportunities Program. The authors are also grateful to M.
Frigo, co-author of FFTW with S.~G.~Johnson \cite{FFTW05}, for
many helpful discussions.

\bibliographystyle{ieeetr}
\bibliography{dct}

\begin{thebibliography}{10}

\bibitem{Johnson07}
S.~G. Johnson and M.~Frigo, ``A modified split-radix {FFT} with fewer
  arithmetic operations,'' {\em {IEEE} Trans. Signal Processing}, vol.~55,
  no.~1, pp.~111--119, 2007.

\bibitem{Kamar89}
I.~Kamar and Y.~Elcherif, ``Conjugate pair fast {Fourier} transform,'' {\em
  Electron. Lett.}, vol.~25, no.~5, pp.~324--325, 1989.

\bibitem{Gopinath89}
R.~A. Gopinath, ``Comment: Conjugate pair fast {Fourier} transform,'' {\em
  Electron. Lett.}, vol.~25, no.~16, p.~1084, 1989.

\bibitem{Qian90}
H.-S. Qian and Z.-J. Zhao, ``Comment: Conjugate pair fast {Fourier}
  transform,'' {\em Electron. Lett.}, vol.~26, no.~8, pp.~541--542, 1990.

\bibitem{Krot92}
A.~M. Krot and H.~B. Minervina, ``Comment: Conjugate pair fast {Fourier}
  transform,'' {\em Electron. Lett.}, vol.~28, no.~12, pp.~1143--1144, 1992.

\bibitem{Vetterli84}
M.~Vetterli and H.~J. Nussbaumer, ``Simple {FFT} and {DCT} algorithms with
  reduced number of operations,'' {\em Signal Processing}, vol.~6, no.~4,
  pp.~267--278, 1984.

\bibitem{Duhamel86}
P.~Duhamel, ``Implementation of ``split-radix'' {FFT} algorithms for complex,
  real, and real-symmetric data,'' {\em {IEEE} Trans. Acoust., Speech, Signal
  Processing}, vol.~34, no.~2, pp.~285--295, 1986.

\bibitem{DuhamelVe90}
P.~Duhamel and M.~Vetterli, ``Fast {Fourier} transforms: a tutorial review and
  a state of the art,'' {\em Signal Processing}, vol.~19, pp.~259--299, April
  1990.

\bibitem{VuducDe00}
R.~Vuduc and J.~Demmel, ``Code generators for automatic tuning of numerical
  kernels: experiences with {FFTW},'' in {\em Proc.~Semantics, Application, and
  Implementation of Code Generators Workshop}, (Montreal), Sept. 2000.

\bibitem{FFTW05}
M.~Frigo and S.~G. Johnson, ``The design and implementation of {FFTW3},'' {\em
  Proc. IEEE}, vol.~93, no.~2, pp.~216--231, 2005.

\bibitem{Wang82}
Z.~Wang, ``A fast algorithm for the discrete sine transform implemented by the
  fast cosine transform,'' {\em {IEEE} Trans. Acoust., Speech, Signal
  Processing}, vol.~30, no.~5, pp.~814--815, 1982.

\bibitem{Chan90}
S.~C. Chan and K.~L. Ho, ``Direct methods for computing discrete sinusoidal
  transforms,'' {\em IEE Proceedings F}, vol.~137, no.~6, pp.~433--442, 1990.

\bibitem{Lee94}
P.~Lee and F.-Y. Huang, ``Restructured recursive {DCT} and {DST} algorithms,''
  {\em {IEEE} Trans. Signal Processing}, vol.~42, no.~7, pp.~1600--1609, 1994.

\bibitem{Yavne68}
R.~Yavne, ``An economical method for calculating the discrete {Fourier}
  transform,'' in {\em Proc.~AFIPS Fall Joint Computer Conf.}, vol.~33,
  pp.~115--125, 1968.

\bibitem{Duhamel84}
P.~Duhamel and H.~Hollmann, ``Split-radix {FFT} algorithm,'' {\em Electron.
  Lett.}, vol.~20, no.~1, pp.~14--16, 1984.

\bibitem{Martens84}
J.~B. Martens, ``Recursive cyclotomic factorization---a new algorithm for
  calculating the discrete {Fourier} transform,'' {\em {IEEE} Trans. Acoust.,
  Speech, Signal Processing}, vol.~32, no.~4, pp.~750--761, 1984.

\bibitem{Lundy07}
T.~Lundy and J.~Van~Buskirk, ``A new matrix approach to real {FFT}s and
  convolutions of length $2^k$,'' {\em Computing}, vol.~80, no.~1, pp.~23--45,
  2007.

\bibitem{Lee84}
B.~G. Lee, ``A new algorithm to compute the discrete cosine transform,'' {\em
  {IEEE} Trans. Acoust., Speech, Signal Processing}, vol.~32, no.~6,
  pp.~1243--1245, 1984.

\bibitem{Wang85}
Z.~Wang, ``On computing the discrete {Fourier} and cosine transforms,'' {\em
  {IEEE} Trans. Acoust., Speech, Signal Processing}, vol.~33, no.~4,
  pp.~1341--1344, 1985.

\bibitem{Suehiro86}
N.~Suehiro and M.~Hatori, ``Fast algorithms for the {DFT} and other sinusoidal
  transforms,'' {\em {IEEE} Trans. Acoust., Speech, Signal Processing},
  vol.~34, no.~3, pp.~642--644, 1986.

\bibitem{Hou87}
H.~S. Hou, ``A fast algorithm for computing the discrete cosine transform,''
  {\em {IEEE} Trans. Acoust., Speech, Signal Processing}, vol.~35, no.~10,
  pp.~1455--1461, 1987.

\bibitem{Arguello90}
F.~Arguello and E.~L. Zapata, ``Fast cosine transform based on the successive
  doubling method,'' {\em Electronic letters}, vol.~26, no.~19, pp.~1616--1618,
  1990.

\bibitem{Kok97}
C.~W. Kok, ``Fast algorithm for computing discrete cosine transform,'' {\em
  {IEEE} Trans. Signal Processing}, vol.~45, no.~3, pp.~757--760, 1997.

\bibitem{Takala00}
J.~Takala, D.~Akopian, J.~Astola, and J.~Saarinen, ``Constant geometry
  algorithm for discrete cosine transform,'' {\em {IEEE} Trans. Signal
  Processing}, vol.~48, no.~6, pp.~1840--1843, 2000.

\bibitem{Puschel03}
M.~P{\"{u}}schel and J.~M.~F. Moura, ``The algebraic approach to the discrete
  cosine and sine transforms and their fast algorithms,'' {\em SIAM
  J.~Computing}, vol.~32, no.~5, pp.~1280--1316, 2003.

\bibitem{Puschel03:dct}
M.~P{\"{u}}schel, ``{Cooley}--{Tukey} {FFT} like algorithms for the {DCT},'' in
  {\em Proc.~IEEE Int'l Conf.~Acoustics, Speech, and Signal Processing},
  vol.~2, (Hong Kong), pp.~501--504, April 2003.

\bibitem{Ahmed74}
N.~Ahmed, T.~Natarajan, and K.~R. Rao, ``Discrete cosine transform,'' {\em IEEE
  Trans. Comput.}, vol.~23, no.~1, pp.~90--93, 1974.

\bibitem{Haralick76}
R.~M. Haralick, ``A storage efficient way to implement the discrete cosine
  transform,'' {\em IEEE Trans. Comput.}, vol.~25, pp.~764--765, 1976.

\bibitem{Chen77}
W.-H. Chen, C.~H. Smith, and S.~C. Fralick, ``A fast computational algorithm
  for the discrete cosine transform,'' {\em IEEE Trans. Commun.}, vol.~25,
  no.~9, pp.~1004--1009, 1977.

\bibitem{Narasimha78}
M.~J. Narasimha and A.~M. Peterson, ``On the computation of the discrete cosine
  transform,'' {\em IEEE Trans. Commun.}, vol.~26, no.~6, pp.~934--936, 1978.

\bibitem{Tseng78}
B.~D. Tseng and W.~C. Miller, ``On computing the discrete cosine transform,''
  {\em IEEE Trans. Comput.}, vol.~27, no.~10, pp.~966--968, 1978.

\bibitem{Yaroslavskii79}
L.~P. Yaroslavskii, ``Shifted discrete {Fourier} transforms,'' {\em Problems of
  Information Transmission}, vol.~15, no.~4, pp.~324--327, 1979.

\bibitem{Makhoul80}
J.~Makhoul, ``A fast cosine transform in one and two dimensions,'' {\em {IEEE}
  Trans. Acoust., Speech, Signal Processing}, vol.~28, no.~1, pp.~27--34, 1980.

\bibitem{Malvar87}
H.~S. Malvar, ``Fast computation of the discrete cosine transform and the
  discrete {Hartley} transform,'' {\em {IEEE} Trans. Acoust., Speech, Signal
  Processing}, vol.~35, no.~10, pp.~1484--1485, 1987.

\bibitem{Li91}
W.~Li, ``A new algorithm to compute the {DCT} and its inverse,'' {\em {IEEE}
  Trans. Signal Processing}, vol.~39, no.~6, pp.~1305--1313, 1991.

\bibitem{Steidl91}
G.~Steidl and M.~Tasche, ``A polynomial approach to fast algorithms for
  discrete {Fourier}-cosine and {Fourier}-sine transforms,'' {\em Math. Comp.},
  vol.~56, no.~193, pp.~281--296, 1991.

\bibitem{Feig92}
E.~Feig and S.~Winograd, ``On the multiplicative complexity of discrete cosine
  transforms,'' {\em IEEE Trans. Info. Theory}, vol.~38, no.~4, pp.~1387--1391,
  1992.

\bibitem{Astola99}
J.~Astola and D.~Akopian, ``Architecture-oriented regular algorithms for
  discrete sine and cosine transforms,'' {\em {IEEE} Trans. Signal Processing},
  vol.~47, no.~4, pp.~1109--1124, 1999.

\bibitem{Guo01}
Z.~Guo, B.~Shi, and N.~Wang, ``Two new algorithms based on product system for
  discrete cosine transform,'' {\em Signal Processing}, vol.~81,
  pp.~1899--1908, 2001.

\bibitem{Plonka05}
G.~Plonka and M.~Tasche, ``Fast and numerically stable algorithms for discrete
  cosine transforms,'' {\em Lin. Algebra and its Appl.}, vol.~394,
  pp.~309--345, 2005.

\bibitem{Frigo99}
M.~Frigo, ``A fast {Fourier} transform compiler,'' in {\em Proc.~ACM SIGPLAN'99
  Conference on Programming Language Design and Implementation (PLDI)},
  vol.~34, (Atlanta, Georgia), pp.~169--180, ACM, May 1999.

\bibitem{Arai88}
Y.~Arai, T.~Agui, and M.~Nakajima, ``A fast {DCT-SQ} scheme for images,'' {\em
  Trans. IEICE}, vol.~71, no.~11, pp.~1095--1097, 1988.

\bibitem{Pennebaker93}
W.~B. Pennebaker and J.~L. Mitchell, {\em {JPEG} Still Image Data Compression
  Standard}.
\newblock New York: Van Nostrand Reinhold, 1993.

\bibitem{Swarztrauber82}
P.~N. Swarztrauber, ``Vectorizing the {FFT}s,'' in {\em Parallel Computations}
  (G.~Rodrigue, ed.), pp.~52--83, New York: Academic Press, 1982.

\bibitem{PressFlaTeu92}
W.~H. Press, B.~P. Flannery, S.~A. Teukolsky, and W.~T. Vetterling, {\em
  Numerical Recipes in {C}: The Art of Scientific Computing}.
\newblock New York, NY: Cambridge Univ.~Press, 2nd~ed., 1992.

\bibitem{VanLoan92}
C.~van Loan, {\em Computational Frameworks for the Fast Fourier Transform}.
\newblock Philadelphia: SIAM, 1992.

\bibitem{RaoYip90}
K.~R. Rao and P.~Yip, {\em Discrete Cosine Transform: Algorithms, Advantages,
  Applications}.
\newblock Boston, MA: Academic Press, 1990.

\bibitem{SorensenJo87}
H.~V. Sorensen, D.~L. Jones, M.~T. Heideman, and C.~S. Burrus, ``Real-valued
  fast {Fourier} transform algorithms,'' {\em {IEEE} Trans. Acoust., Speech,
  Signal Processing}, vol.~35, no.~6, pp.~849--863, 1987.

\bibitem{GenSan66}
W.~M. Gentleman and G.~Sande, ``Fast {F}ourier transforms---for fun and
  profit,'' {\em Proc.~AFIPS}, vol.~29, pp.~563--578, 1966.

\bibitem{Schatzman96}
J.~C. Schatzman, ``Accuracy of the discrete {F}ourier transform and the fast
  {F}ourier transform,'' {\em SIAM J.~Scientific Computing}, vol.~17, no.~5,
  pp.~1150--1166, 1996.

\bibitem{Tasche00}
M.~Tasche and H.~Zeuner, {\em Handbook of Analytic-Computational Methods in
  Applied Mathematics}, ch.~8, pp.~357--406.
\newblock Boca Raton, FL: CRC Press, 2000.

\bibitem{CrochiereOp75}
R.~E. Crochiere and A.~V. Oppenheim, ``Analysis of linear digital networks,''
  {\em Proc. {IEEE}}, vol.~63, no.~4, pp.~581--595, 1975.

\bibitem{ShaoJo08-dct4}
X.~Shao and S.~G. Johnson, ``Type-{IV} {DCT}, {DST}, and {MDCT} algorithms with
  reduced numbers of arithmetic operations,'' {\em Signal Processing}, 2008.
\newblock In press.

\end{thebibliography}

\end{document}